\documentclass[reqno,11pt]{amsart}
\usepackage{amsmath, latexsym, amsfonts, amssymb, amsthm, amscd}

\numberwithin{equation}{section}

\newtheorem{theorem}{Theorem}

\newtheorem{proposition}[theorem]{Proposition}
\newtheorem{corollary}[theorem]{Corollary}
\newtheorem{rem}{Remark}
\newtheorem{definition}[theorem]{Definition}

\DeclareMathSymbol{\leqslant}{\mathalpha}{AMSa}{"36} 
\DeclareMathSymbol{\geqslant}{\mathalpha}{AMSa}{"3E} 
\DeclareMathSymbol{\eset}{\mathalpha}{AMSb}{"3F}     

\newcommand{\sumtwo}[2]{\sum_{\substack{#1 \\ #2}}} 

\DeclareMathOperator{\sign}{sign}

\title{Polymer Pinning at an interface}

\author{Nicolas Petrelis}

\address{Laboratoire de Math\'ematiques Raphael Salem, site Colbert, 76821 Mont Saint Aignan , France
CNRS U.M.R. 6085}

\email{nicolas.petrelis\@@etu.univ-rouen.fr}

\date{\today}

\begin{document}

\begin{abstract}
We consider a model of hydrophobic homopolymer in interaction with an interface between oil and water.
The configurations of the polymer are given by the trajectories of a simple symmetric random walk
$\left(S_{i}\right)_{i\geq 0}$. On the one hand the hydrophobicity of each monomer tends to delocalize
the polymer in
the upper half plane, that is why we define
$h$, a non negative energetic factor that the chain gains for every monomer
in the oil (above the origin). On the other hand the chain receives a random price (or penalty) on
crossing the interface.
At site $i$ this  price is given by
$\beta\left(1+s\zeta_{i}\right)$, where $(\zeta_{i})_{i\geq 1}$ is a sequence of i.i.d. centered
 random variables, and $(s,\beta)$ are two non negative parameters. Since the price is positive on the average,
the interface
 attracts the polymer and a localization effect may arise. We transform the measure
of
 each trajectory with the hamiltonian $\beta\sum_{i=1}^{N} (1+s \zeta_{i})\boldsymbol{1}_{\{S_{i}=0\}}+h
 \sum_{i=1}^{N} \sign(S_{i})$, and study the critical curve $h^{s}_{c}\left(\beta\right)$ that divides
  the phase spaces in a localized and a delocalized area.

It is not difficult to show that $h_{c}^{0}(\beta)\leq h^{s}_{c}\left(\beta\right)$ for every $s\geq 0$,
but in this article we give a method to improve in a quantitative way this lower
bound. To that aim, we transform the strategy developed by Bolthausen and Den Hollander in \cite{BDH} on
taking into account the fact
that the chain can target the sites
where it comes back to the origin. Then we deduce from this last
result a corollary in terms of pure pinning model, namely
with the hamiltonian $\sum_{i=1}^{N} (-u+s \zeta_{i})\boldsymbol{1}_{\{S_{i}=0\}}$
we find a lower bound of the critical curve $u_{c}(s)$ for small $s$.
 In this situation, we improve the existing lower bound of Alexander and Sidoravicius \cite{SidAlex}.

\bigskip


\bigskip
\end{abstract}

\maketitle

\section{Introduction and results}
\subsection{the model}
We consider a simple random walk $(S_{n})_{n\geq 0}$, defined as $S_{0}=0$ and
$S_{n}= \sum_{i=1}^{n}X_{i}$ where
$(X_{i})_{i\geq 1}$  is a sequence of iid bernouilli trials verifying $P(X_{1}=\pm 1)=1/2$.
We denote by $\Lambda_{i}=\sign(S_{i})$ if $S_{i}\neq 0$,\ \ $\Lambda_{i}=\Lambda_{i-1}$ otherwise.
We also define $(\zeta_{i})_{i\geq 1}$ a sequence of iid random variables non a.s. equal to $0$,
 verifying\ \ $\mathbb{E}
\left(e^{\lambda|\zeta_{1}|}\right)<\infty$ for every $\lambda>0$
and $\mathbb{E}\left(\zeta_{1}\right)=0$.

 Now let $h\geq 0$, $s\geq 0$ and for each trajectory of the random walk we define the following hamiltonian
$$H_{N,\beta,h}^{\zeta,s}(S)=\beta \sum_{i=1}^{N}\left(
1+s\zeta_{i}\right) \boldsymbol{1}_{\{S_{i}=0\}}+h\sum_{i=1}^{N}\Lambda_{i}$$
With this hamiltonian we perturb the law of the random walk as follow
$$\frac{dP_{N,\beta,h}^{\zeta,s}}{dP}\left(S \right)=
\frac{\exp\left(H_{N,\beta,h}^{\zeta,s}(S)\right)}{Z_{N,\beta,h}^{\zeta,s}}$$
This new measure $P_{N,\beta,h}^{\zeta,s}$ is called polymer
measure of size $N$. Under this measure two sorts of trajectories
are "a priori" favored. On the one hand the localized
trajectories, often coming back to the origin to receive some
positive pinning rewards along the x axis. On the other hand,
trajectories called delocalized, spending most time in the upper
half plane and both favored by the second term of the hamiltonian
and the fact that they are much more numerous than those staying
close to the origin. So a competition between these two possible
behaviors arises.

\subsection{previous results and physical motivations}
Systems of random walk attracted by a potential at an interface
are closely studied at this moment (see \cite{Mont}, \cite{JRV}). One of the major issue about
that subject consists in understanding better the influence of a
random potential compared to a constant one. Namely, if it seems
intuitively clear that a random potential has a stronger power of
attraction than a constant one of same expectation, it is much
more complicated to quantify this difference.

 In the present work, we consider a potential
at the interface and also the fact that the polymer prefers lying
in the  upper half plan than in the lower one. That type of system has been studied numerically in \cite{JSW},
and
  can describe for example the situation of an hydrophobic
 homopolymer 
 at an interface between oil and
water. Close to this horizontal separation between the two
solvents, some very small droplets of a third solvent
(microemulsions) are put and have a big power of attraction on the monomers
composing our chain. So the pinning prices our chain can receive
when it comes back to the origin represent the attractive
emulsions our polymer can touch close to the interface.

  We expose here precise theoretical results about the critical
 curve arising from this system. We investigate new strategies of
 localization for the polymer consisting in targeting the sites
 where it comes back to the interface, and we find an explicit
 lower bound of the critical curve strictly above the non random
 one. 

 Our result covers, as a limit case as $h$ goes to infinity the wetting transition model.
  Effectively in the last ten years
  the wetting problem, namely the case of
  a polymer interacting with an (impenetrable) interface has attracted a lot of interest
 since it can be regarded as a Polland Sheraga model of the DNA strand  (see \cite{Mont}).
 The localization transition
 with a constant disorder occurs for the pinning reward $\log 2$, and a lot of
 questions arising
 from this first result are linked with the effect of a small random
 perturbation add to the price $\log 2$.  Moreover, with the constant pinning reward $\log 2$
  the simple random walk conditioned to stay positive
 has the same law than the reflected random walk (see \cite{Yosh}). That is why, to study
 the wetting model around the pinning price $\log(2)$, it suffices to consider the pure
 pinning model, namely a reflected random walk pinned at the origin by small random variables.

This last model has therefore been closely studied,
for example in \cite{JRV} a particular type of positive potential has been
considered and a criterium has been given to decide for every
disorder realization if it localizes the polymer or not. But a very
difficult question consists in estimating, for small $s$, the critical delocalization
 average
$u_{c}(s)$
of an iid disorder of type $-u+s\zeta_{i}$ with $\zeta_{i}$ centered of variance $1$ $\big($namely Var$\left(
-u+s\zeta_{i}\right)=s^{2}$ $\big)$. The annealed critical curve is given by
$u_{a}(s)=\log E\left(\exp(s\zeta_{i})\right)\stackrel{s \to 0} {\sim} s^{2}/2$ $\Big($even $=s^{2}/2$ when
$\zeta_{i} \stackrel{D}{=} \emph{N}(0,1)\Big)$ and
verifies as usual $u_{c}(s)\leq u_{a}(s)$. 
  In the last $20$ years there has been a lot of activity on this question, mostly from
the physicists side and it is now widely believed that $u_{c}(s)$ behaves as $s^{2}/2$ but
it is still an open question wether $u_{c}(s)=s^{2}/2$ (see
\cite{For}) for $s$ small or $u_{c}(s)<s^{2}/2$ for every $s$ (see \cite{Der} or \cite{Kal}).

However up to now the only rigorous thing that has been proved is
  in \cite{SidAlex}, where Sidoravicious and Alexander have
 studied a general class of random walk pinned either by an interface between
 two solvents or by an impenetrable wall. If we apply their results in our case it gives
that for iid centered $(\zeta_{i})_{i\geq0}$ of fixed positive variance, the quenched quantity $u_{c}(s)$
 is strictly
larger than the non disordered one $u_{c}(0)$. In this paper, the new localization strategies
 we develop allows us to go further
on giving an upper bound of $u_{c}(s)$ of type $-c s^{2}$, which has the same scale than the
annealed lower bound.

\subsection{the free energy}
To decide for fixed parameters if our system is localized or not
we introduce the free energy called $\Psi^{s}(\beta,h)$ and
defined as
$$\Psi^{s}(\beta,h)=\lim_{N \to \infty}\frac{1}{N} \log Z_{N,\beta,h}^{\zeta,s}$$
This limit $\Psi^{s}(\beta,h)$ is not random any more and occurs
$\mathbb{P}$ almost surely in $\zeta$ and $\mathbb{L}^{1}$. The proof of that
sort of convergence is well known (see \cite{Giac} or \cite{BDH}).
This free energy can easily be bounded from below on computing it
on a restriction of the trajectories set. That way we denote by $D_{N}$
the set $\{S\ :\
 S_{i}>0\  \forall i\in \{1,..N\}\}$.
For each trajectory of $D_{N}$ the hamiltonian is equal to $h N$ since the chain stays in the upper half
plane and never comes back to the origin. Moreover $P\left(D_{N}\right)\sim c/N^{1/2}$ as $N$ goes to $\infty$.
Hence
\begin{equation*}
\Psi^{s}(\beta,h)\geq \liminf_{N \to \infty}\frac{1}{N}\log E\left(e^{hN} \boldsymbol{1}_{\{D_{N}\}}\right)
\geq h+\liminf_{N \to \infty}\frac{\log\left(P\left(D_{N}\right)\right)}{N}
\geq h
\end{equation*}
so the free energy is always larger than $h$, and from now on we
will say that the polymer is delocalized if $\Psi^{s}(\beta,h)=h$,
because the utterly delocalized trajectories of $D_{N}$ give us
the whole free energy whereas it will be delocalized if
$\Psi^{s}(\beta,h)>h$.

 This separation between localized and delocalized regime seems a
 bit raw, because many trajectories come back only a few times to
 the origin and should also be called delocalized since they spend
 almost all their time in the upper half plane. So taking only into
 account the utterly delocalized trajectories could be
not sufficient. But it is in fact because for convexity reasons,
in all the localized phase the chains come back to the origin a
positive density of times. Another result can help us to
understand the localization effect. It is due to Sinai in \cite{Sinai}
 and with the same technics we can control the vertical expansion
 of the chain in the localized area. That way we transform a bit the
 hamiltonian which becomes $\left( \beta \sum_{i=1}^{N}\left(1+s\zeta_{N-i}\right)
 \boldsymbol{1}_{\{S_{i}=0\}}+h\sum_{i=1}^{N}\left(\Lambda_{i}\right)\right)$ so that the disorder is fixed
 in the neighborhood of $S_{N}$. Notice that the free energy is not modified by this transformation
 and allows us to say for $\Psi^{s}(\beta,h)>0$, $\epsilon>0$ and every realization of the
 disorder $\zeta$ that
there exists a constant $C_{\zeta}^{\epsilon}>0$, $\mathbb{P}$ almost surely finite
verifying for every $L\geq 0$ and $N\geq 0$
$$P_{N,\beta,h}^{\zeta,s}\left(|S_{N}|>L\right)\leq C_{\zeta}^{\epsilon}
\exp\left(-\left(\Psi^{s}(\beta,h)-\epsilon\right)L\right)$$
 This result can not occur if we keep the original hamiltonian
 because the disorder is not fixed close to $S_{N}$. As a consequence
 we meet almost surely arbitrary long stretches of negative rewards
that push rarely but sometimes $S_{N}$ far away from the interface. 

Some pathwise results have also been proved in the delocalized area for polymer systems. In our case
we can use the method developed in the last part of \cite{BiskDH} to prove that $\mathbb{P}$ almost
surely in $\zeta$ and for every $K>0$,
$\lim_{N \to \infty}E_{N,\beta,h}^{\zeta,s}\left(\sharp\{i\in\{1,..,N\} : S_{i}>K\}/N\right)=1$. These results
allow us to understand more deeply what localization and delocalization mean.


Now we want to transform the hamiltonian, in order to simplify the localization condition. In that way
notice that
$$\Psi^{s}(\beta,h)-h=\lim_{N\to \infty}\frac{1}{N}\log\left(E\left( \exp\left(\beta \sum_{i=1}^{N}\left(
1+s\zeta_{i}\right) \boldsymbol{1}_{\{S_{i}=0\}}+h\sum_{i=1}^{N}\left(\Lambda_{i}-1\right)\right)\right)\right)$$
so we put $\Phi^{s}(\beta,h)=\Psi^{s}(\beta,h)-h$, the delocalization condition becomes $\Phi^{s}(\beta,h)=0$
 and the
localization one $\Phi^{s}(\beta,h)>0$. To finish with these new notations we denote $\Delta_{i}=1$
if $\Lambda_{i}=-1$
and $\Delta_{i}=0$ if $\Lambda_{i}=1$. The hamiltonian becomes
$$H_{N,\beta,h}^{\zeta, s}\left(S\right)=\beta \sum_{i=1}^{N}\left(
1+s\zeta_{i}\right) \boldsymbol{1}_{\{S_{i}=0\}}-2h\sum_{i=1}^{N}\Delta_{i}$$
and we keep
$Z_{N,\beta,h}^{\zeta,s}=E\left(e^{H_{N,\beta,h}^{\zeta,s}}\right)$,
so we have
$$\Phi^{s}(\beta,h)=\lim_{N \to \infty}\frac{1}{N} \log Z_{N,\beta,h}^{\zeta,s}$$
This function $\Phi^{s}$ is convex and continue in both variables, non decreasing in $\beta$ and non increasing in
$h$.
In this paper we are particularly interested in the critical curve of the system, namely the curve that
divides the phases space ($h$,$\beta$) in a delocalized zone, and a localized one. But before defining this
curve precisely, it is helpful to consider the non disordered case ($s=0$) , which is much simpler to perform and
provides intuitions about what happens in the disordered case ($s\neq 0$).
\subsection{the critical curve}
Above the critical curve the system will be delocalized and localized below, in appendix C)
we compute the equation of this curve when $s=0$, we obtain
\begin{align}\label{pro:nonrand}
 \nonumber h^{0}_{c}:\ [0,&\log(2))\rightarrow \mathbb{R}\\
 \beta &\longrightarrow h_{c}^{0}\left(\beta\right)=\frac{1}{4}\log\left(1-4\left(1-e^{-\beta}\right)^{2}\right)
\end{align}
So the curve is increasing, convex and goes to $\infty$ when $\beta$ goes to $\log(2)$ from the left.
But when $\beta\geq \log(2)$ the system is always localized, in fact as large as $h$ is chosen the free energy
 remains strictly positive, that is why this critical curve is only defined on $[0,\log(2))$ (see Fig $1$).

This lets us think that when $s\neq 0$, the critical curve should
have a form of the same type as \eqref{pro:nonrand}. Notice also
that $h_{c}^{0}\left(\beta\right)\sim \beta^{2}$as $\beta$ goes to $0$.

\begin{proposition}\label{pro:sec}
For every $s \geq 0$ and $\beta \geq 0$ there exists
$h^{s}_{c}(\beta) \in [0,+\infty]$ 
 such that for every $h< h^{s}_{c}(\beta)$ the free energy
 $\Phi^{s}(\beta,h)$ is strictly positive, whereas $\Phi^{s}(\beta,h)=0$ if $h \geq
 h^{s}_{c}(\beta)$ . 
This function $h_{c}^{s}(\beta)$ is convex, increasing in
$\beta$, hence for every $s\geq 0$ there exists $\beta_{0}(s)\in [0,\infty]$
verifying $h^{s}_{c}(\beta)< +\infty$ if $\beta<\beta_{0}(s)$
and $h^{s}_{c}(\beta)= +\infty$ if $\beta>\beta_{0}(s)$.
\end{proposition}
We will prove also that for every $s\geq 0$ the non
disordered critical curve $h_{c}^{0}(\beta)$ is a lower bound of
$h_{c}^{s}(\beta)$. As a consequence $\beta_{0}(s)\leq
\beta_{0}(0)=\log(2)$
\begin{rem}\label{remarque sur courbe}
The case $\beta=\beta_{0}(s)$ remains open, more precisely two
different behavior of the curve may occur. Either $\lim_{\beta \to
\beta_{0}^{-}(s)}h_{c}\left(\beta\right)=+\infty$, or there exists $h_{0}^{s}<\infty$ such that $\lim_{\beta
\to \beta_{0}^{-}(s)}h_{c}\left(\beta\right)=h_{0}^{s}$ and by
continuity of $\Phi^{s}$ in $\beta$ we have
$\Phi(\beta_{0}(s),h_{0}^{s})=0$ and
$h_{c}\left(\beta_{0}(s)\right)=h_{0}^{s}$.
\end{rem}

We Find an upper bound of $h_{c}^{s}(\beta)$ as usual, on
computing the annealed free energy, which is by the Jensen
inequality an upper bound of the quenched free energy. The
annealed system gives birth to a critical curve
($h^{s}_{ann.c}(\beta)$) which is an upper bound of the quenched
critical curve. 
The annealed free energy is
given by
$$\Phi^{s}_{ann.}(h,\beta)=\lim_{N\rightarrow\infty}\frac{1}{N}\log\,E\mathbb{E}\left(\exp\left(\beta\sum_{i=1}^{N}
\left(1+s\zeta_{i}\right) \boldsymbol{1}_{\{S_{i}=0\}}-2h
\sum_{i=1}^{N}\Delta_{i}\right)\right)$$ hence if we integrate
over $\mathbb{P}$ we obtain
$$\Phi^{s}_{ann.}(h,\beta)=\lim_{N\rightarrow\infty}\frac{1}{N}\log\,E\left(\exp\left(
\left(\beta+\log\mathbb{E}(e^{\beta s
\zeta_{1}})\right)\sum_{i=1}^{N} \boldsymbol{1}_{\{S_{i}=0\}}-2h
\sum_{i=1}^{N}\Delta_{i}\right)\right)$$
 Finally $\Phi^{s}_{ann.}(h,\beta)=\Phi^{0}(h,\beta+\log\mathbb{E}(e^{\beta s
 \zeta_{1}}))$ and the annealed critical curve can be expressed with the help of
the non disordered one, namely if we call
 $\beta_{ann}^{s}$ the only solution of\  $\beta+\log\mathbb{E}(e^{\beta s \zeta_{1}})=\log
 2$,
 for every $\beta \in [0,\beta_{ann}^{s})$ the value of the annealed
critical curve is
$h_{ann.c}^{s}(\beta)=h_{c}^{0}\left(\beta+\log\mathbb{E}\left(e^{\beta
s \zeta_{1}}\right)\right)$ (see Fig $1$).

 Once again notice that the annealed critical
curve verifies $h_{ann.c}^{s}(\beta)\sim \beta^{2}$as $\beta$ goes
to $0$.

\subsection{The disordered model}
Here comes the main part of the paper, we develop a new strategy
to find a lower bound on the quenched critical curve. A strategy
to find that kind of lower bound consists in computing the free
energy on a particular restriction of the trajectories, namely, in
the localized area, trajectories that often come back to the
origin (\cite{BodGiac}). Here we are going to develop another
method, that consists in transforming (using radon Nikodym
densities) the law of the excursions out of the origin. First (as
done
 in \cite{BDH}) we constrain the chain to come back to the origin a positive density
of times, but without targeting the sites of the $x$ axes it will
touch. Then we make the chain choose at each excursion a
trajectory law adapted to the local environment.

 Notice first that proposition 1  tells us that for every $s\geq 0$
 and $\beta \geq \log(2)$ we have $h_{c}^{s}(\beta)=\infty$ hence
 in any case the critical curve is not defined after $\log{2}$,
 that is why, from now on we only consider the case $\beta\leq \log(2)$.
\begin{theorem}\label{theo}
If $Var(\zeta_{1})\in (0, \infty)$, there exists two strictly positive constant $c_{1}$ and
$c_{2}$ such that for every $s\leq c_{1}$ and
$\beta \in [0,\log{2}-c_{2} s^{2} \beta^{2})$,
 we can bound from below the critical
curve as follow
$$h_{c}^{s}(\beta)\geq -\frac{1}{4}\log\left(1-4\left(1-e^{-\beta-c_{2} s^{2} \beta^{2}}\right)^{2}\right)=m^{s}(\beta)$$
\end{theorem}
\begin{rem}
This lower bound is strictly above the non disordered one (see
proposition 1 and Fig $1$) when $s>0$.
\end{rem}


\vspace{2cm}

\setlength{\unitlength}{0.35cm}

\begin{picture}(30,15)(-7,0)

  \put(0,0){\vector(1,0){28}}
  \put(0,0){\vector(0,1){15}}
 \put(25,0){\line(0,15){15}}
  \put(-.8,-1.3){$0$}
  \put(29,-.3){$\beta$}
  \put(25,-1.3){$\log 2$}
  \put(-1,16){$h$}
  \put(14,-1.3){$\beta_{ann}^{s}$}

  \put(16,3.2){\line(0,12){9.9}}
   \put(17,3.9){\line(0,12){9.2}}
   \put(18,4.7){\line(0,12){8.4}}
   \put(19,5.6){\line(0,12){7.5}}
   \put(20,6.9){\line(0,12){6.2}}
   \put(21,8.7){\line(0,12){4.4}}

   \put(14,2.3){\line(0,12){6.3}}
   \put(13,1.9){\line(0,12){4.2}}
   \put(12,1.5){\line(0,12){3.2}}
   \put(11,1.2){\line(0,12){2.3}}
   \put(10,1){\line(0,12){1.7}}
   \put(9,0.8){\line(0,12){1.3}}
   \put(8,0.6){\line(0,12){0.9}}
   \put(7,0.4){\line(0,12){0.7}}
   \put(6,0.3){\line(0,12){0.5}}

  {\thicklines
   \qbezier[90](0,0)(24,0)(24.6,13)
  }
  {\thicklines
   \qbezier(15,0)(15,5)(15,15)
  }
{\thicklines
   \qbezier(22.3,0)(22.3,5)(22.3,15)
  }
  {\thicklines
   \qbezier(25,0)(25,5)(25,15)
  }
{\thicklines
   \qbezier[30](0,0)(14,0)(14.6,13)
}
{\thicklines
   \qbezier(0,0)(21,0)(22,13)
  }

  \put(0,-4.5){\small
             Fig.\ 1:{\thicklines
   \qbezier[6](3,0)(5,0)(7,0)}

   {\thicklines
   \qbezier(3,-3)(5,-3)(7,-3)}
{\thicklines
   \qbezier[12](3,-6)(5,-6)(7,-6)}}
\put(12,-4.5){\small
             $h_{ann.c}^{s}(\beta)$}
\put(12,-10.5){\small
             $h_{c}^{0}(\beta)$}
\put(12,-7.5){\small
             $m^{s}(\beta)$}
\put(20.5,-6.5){\line(10,0){3}}
\put(20.5,-8.5){\line(10,0){3}}
\put(20.5,-6.5){\line(0,-3){2}}
\put(21.5,-6.5){\line(0,-3){2}}
\put(22.5,-6.5){\line(0,-3){2}}
\put(23.5,-6.5){\line(0,-3){2}}
\put(25,-7.7){\small
             possible location of $h_{c}^{s}(\beta)$}
\put(15,0){\circle*{.45}}
\put(25,0){\circle*{.45}}
\put(22.3,0){\circle*{.45}}
\end{picture}

\vspace{5cm}

\begin{rem}
 The possible values of $c_{1}$ and $c_{2}$ depend on the law of $\zeta_{1}$. For example, as showed in the proof
 if $\mathbb{P}\left(\zeta_{1}
 >0\right)=1/2$ and $\mathbb{E}\left(\zeta_{1} \boldsymbol{1}_{\{\zeta_{1}>0\}}\right)=1$, the values
 $c_{1}=1$ and $c_{2}=1/(5\times 2^{14})$ suit.
   But with other conditions the strategy to obtain the lower bound remains the same.
 \end{rem}
\begin{rem}
The precise value of $c_{2}$ $\left(1/\left(5 \times 2^{14}\right)\right)$ could certainly be improved,
on building more complicated law of return
to the origin. For example on building a law of return
to the origin that depends more deeply on the environment (taking into account $\zeta_{i+2}, \zeta_{i+4}$ etc...).
The computations would be quite more complicated and our aim here is not to optimize the value of $c$ but
to expose a simple strategy that improves the non
disordered lower bound of a term $c s^{2} \beta^{2}$ with $c>0$.
\end{rem}

\subsection{The pure pinning model}
The pure pinning model is a bit different from our previous one, the $h$ term of entropic repulsion vanishes
and we consider pinning rewards at the origin of the form $-u+s\zeta_{i}$ with
$u \geq 0$. The corresponding
hamiltonian is
$$H_{N,s}^{\zeta,u }\left((i,S_{i})_{i\in \{0,..,N\}}\right)= \sum_{i=1}^{N}\left(
-u+s\zeta_{i}\right) \boldsymbol{1}_{\{S_{i}=0\}}$$ In that case,
the condition of localization and delocalization in term of free
energy remains the same and we have a critical $u$ called
$u_{c}(s)$ such that for $u\geq u_{c}(s)$ the system is
delocalized, whereas for $u<u_{c}(s)$ it is localized. Recall also that if
$Var(\zeta_{1})=1$ the annealed case tells us that $u_{c}(s)\leq u_{c}^{ann}(s)\sim_{s\to 0} s^{2}/2$.
 Now, a corollary of our theorem gives us a lower bound on
$u_{c}(s)$ which has the good scale.
\begin{corollary}\label{cor}
If $Var(\zeta_{1})\in (0, \infty)$, there exists two strictly positive constant $c_{3}$ and
$c_{4}$ such that for every $s\leq c_{3}$
 $$u_{c}(s)\geq c_{4}\  s^{2}$$
 \end{corollary}
\begin{rem}
Once again the values of $c_{3}$ and $c_{4}$ depend on the law of $\zeta_{1}$. We will keep in our proof the
conditions of remark $3$ concerning $\zeta_{1}$. The values $c_{3}=\log 2$ and $c_{4}=1/(5 \times 2^{16})$ suit.
\end{rem}

\section{Proof of theorem and proposition}

\subsection{\textit{Proof of Proposition \ref{pro:sec}}}
First define for every $\beta\geq 0$ and $s\geq 0$ the set $J_{\beta}^{s}=\{h\geq 0
\ \ \text{such that}\ \ \Phi^{s}\left(\beta,h\right)= 0\}$
. We put $h_{c}^{s}(\beta)$ the
lower bound of $J_{\beta}^{s}$.
Then recall that $\Phi$ is continuous, not increasing in $h$, and positive hence the set
$J_{\beta}^{s}$ can be written $[h_{c}^{s}\left(\beta  \right), +\infty)$ (when it is not empty). Moreover $\Phi$ is
not decreasing in $\beta$ because $\Phi^{s}(0,h)=0$ for every $h\geq 0$, $\Phi\left(\beta,h\right)\geq 0$
for every $\beta$ and $\Phi$ is convex in $\beta$.
 So if $\beta_{1}\geq \beta_{2}$ we have $J_{\beta_{1}}^{s} \subset J_{\beta_{2}}^{s}$. This
gives us the fact that $h_{c}(\beta)$ is not decreasing, and we
put $\beta_{0}(s)=\sup\{\beta \geq 0 : J_{\beta}^{s}\neq \emptyset\}$.
The annealed computation shows us that
$\beta_{0}(s)>0$ because $\Phi^{s}(h,\beta) \leq
\Phi^{s}_{ann.}(h,\beta)$. Thus $J_{ann.\beta} \subset J_{\beta}$
and $\beta_{0}(s) \geq \beta_{ann}^{s}>0$.
 Now we want to prove that $h_{c}(\beta)$ is convex. That way it is continuous
on the interval $[0, \beta_{0}(s))$.\\ 
To prove this convexity we put $0<a<b$ and $\lambda \in [0,1]$. So remark that
$$H_{N,\ \lambda a +\, (1-\lambda) b,\ \lambda h_{c}^{s}(a)+\, (1-\lambda)h_{c}^{s}(b)}^{\zeta,s}=
H_{N,\ \lambda a,\ \lambda h_{c}^{s}(a)}^{\zeta,s}\  +H_{N,\ (1-\lambda) b,\ (1-\lambda) h_{c}^{s}(b)}^{\zeta,s}$$
hence by holder inequality
\begin{align*}
\frac{1}{N}\log E \left(\exp\left(Z_{N,\ \lambda (a,h_{c}^{s}(a)) +\,(1-\lambda) (b,h_{c}^{s}(b))}^{\zeta,s}\right)
\right)
&\leq \frac{\lambda}{N}\log E \left(\exp\left(Z_{N,\ a,\ h_{c}^{s}(a)}^{\zeta,s}\right)\right)\\
&+ \frac{1-\lambda}{N}\log E \left(\exp\left(Z_{N,\ b,\ h_{c}^{s}(b)}^{\zeta,s}\right)\right)
\end{align*}
so as $N$ goes to infinity the two terms of the rhs goes to zero
because by continuity of $\Phi$ in $h$ we have
$\Phi(a,h_{c}^{s}(a))=\Phi(b,h_{c}^{s}(b))=0$. Hence
$\Phi^{s}(\lambda a +(1-\lambda) b, \lambda
h_{c}^{s}(a)+(1-\lambda)h_{c}^{s}(b))=0$ and $h_{c}^{s}(\lambda a
+(1-\lambda) b) \leq \lambda
h_{c}^{s}(a)+(1-\lambda)h_{c}^{s}(b)$. This completes the proof.

 Now it remains to give a short proof of the fact that
 $h_{c}^{s}(\beta)\geq h_{c}^{0}\left(\beta\right)$ for every $s\geq 0$. We will in fact
 prove that for $s\geq 0$, $\beta\geq 0$ and $h<
 h_{c}^{0}\left(\beta\right)$ the free energy $\Phi^{s}(\beta,h)>0$. This
 will be sufficient to complete the proof. Hence notice that for
 fixed $(\beta,h)$ the function $\Phi^{s}(\beta,h)$ is convex in
 $s$ since it is the limit as $N$ goes to infinity of the function sequence
  $\Phi_{N}^{s}(\beta,h)=\mathbb{E}\left(1/N\log E\left((\exp\left(H_{N,\beta,h}^{\zeta,s}\right)\right)\right)$
which are convex in s. Moreover for every $N>0$,
$\Phi_{N}^{s}(\beta,h)$ can be derived in s and this gives
\begin{equation*}
\frac{\partial \Phi_{N}^{s}(\beta,h)}{\partial
s}=\frac{1}{N}\mathbb{E}\left(\frac{E\left(\beta
\sum_{i=1}^{N}\zeta_{i} \boldsymbol{1}_{\{S_{i}=0\}}
\exp\left(H_{N,\beta,h}^{\zeta,s}\right)\right)}{E\left(
\exp\left(H_{N,\beta,h}^{\zeta,s}\right)\right)}\right)
\end{equation*}
But when $s=0$ the hamiltonian does not depend on the disorder
($\zeta$) any more, so by Fubini Tonelli and the fact that the
$\zeta_{i}$ are centered we can write
\begin{equation*}
\frac{\partial \Phi_{N}^{s}(\beta,h)}{\partial
s}\Bigg|_{s=0}=\frac{1}{N}\frac{E\left(\beta
\sum_{i=1}^{N}\mathbb{E}\left(\zeta_{i}\right)
\boldsymbol{1}_{\{S_{i}=0\}}
\exp\left(H_{N,\beta,h}^{\zeta,0}\right)\right)}{E\left(
\exp\left(H_{N,\beta,h}^{\zeta,0}\right)\right)}=0
\end{equation*}
hence the convergence of $\Phi_{N}$ to $\Phi$ and their convexity
allow us to say
\begin{equation*}
\frac{\partial_{right} \Phi^{s}\left(\beta,h\right)}{\partial
s}\Bigg|_{s=0}\geq \lim_{N \to \infty}\frac{\partial_{right}
\Phi_{N}^{0}\left(\beta,h\right)}{\partial s}\Bigg|_{s=0}=0
\end{equation*}
so since $\Phi^{s}\left(\beta,h\right)$ is convex in $s$ we can
conclude that it is not decreasing on $[0,\infty)$. Hence for
every $s\geq 0$, $\Phi^{s}\left(\beta,h\right)\geq
\Phi^{0}\left(\beta,h\right)>0$. That is why $h_{c}^{s}\left(\beta\right)\geq h_{c}^{0}\left(\beta\right)$.

To finish with this proof, we show that $h_{c}^{s}\left(\beta\right)$ is increasing in $\beta$.
In fact since $h_{c}^{s}\left(0\right)=0$ and
$h_{c}^{s}\left(\beta\right)\geq h_{c}^{0}\left(\beta\right)>0$ for $\beta>0$ the convexity of
$h_{c}^{s}\left(\beta\right)$ gives us the result. \qed

\medskip
\noindent
\subsection{\textit{Proof of Theorem \ref{theo}}}
In the following we consider $h>0$,\ $\beta\leq\log(2)$,\ $\mathbb{P}(\zeta_{1}>0)=1/2$,
\ $\mathbb{E}\left(\zeta_{1}
\boldsymbol{1}_{\{\zeta_{1}>0\}}\right)=1$ and $s\leq1$.
\subsubsection{STEP 1: transformation of the excursions law}
\begin{definition}
From now on we will call $i_{j}$ the site where the $j^{th}$ return to the origin
takes place, so $i_{0}=0$ and $i_{j}=\inf\{i>i_{j-1}:S_{i}=0\}$ and $\tau_{j}=i_{j}-i_{j-1}$ is
 the length of the $j^{th}$ excursion out of the origin. We also call $l_{N}$ the number of return to the
  origin before time $N$.
\end{definition}
Thus by independence of the excursions signs we can rewrite the partition function as


\begin{equation}\label{HN}
H_{N}=E\left(\exp\left(\beta s\sum_{j=1}^{l_{N}}
\zeta_{i_{j}}\right)\exp\left(\beta l_{N}\right)
\prod_{j=1}^{l_{N}}\left(\frac{1+\exp\left(-2h\tau_{j}\right)}{2}\right)
\left(\frac{1+\exp\left(-2h\left(N-i_{l_{N}}\right)\right)}{2}\right)\right)
\end{equation}
Now we want to transform the law of excursions out of the origin
to constrain the chain to come back to zero a positive density of
times. That way we introduce $P_{\alpha,h}^{\beta}$ the law of an homogeneous positive recurrent markov
process, whose excursion law are given by
\medskip
\begin{equation}\label{mes}
\forall n \in \mathbb{N}-\{0\}\ \ \ \ \
P_{\alpha,h}^{\beta}\left(\tau_{1}=2n\right)=\left(\frac{1+\exp\left(-4hn\right)}{2}\right)\alpha^{2n}
\frac{P\left(\tau=2n\right)}{H_{\alpha,h}^{\beta}}\exp\left(\beta\right)
\end{equation}
where $H_{\alpha,h}^{\beta}$ can be computed as follow
\begin{equation}\label{H}
H_{\alpha,h}^{\beta}=\sum_{i=1}^{\infty}\frac{\exp\left(-4hi\right)+1}{2}e^{\beta}\alpha^{2i}
P\left(\tau=2i\right)
=e^{\beta}\left(1-\frac{\sqrt{1-\alpha^{2}}+\sqrt{1-e^{-4h}\alpha^{2}}}{2}\right)
\end{equation}
\medskip
Notice also that the function we are considering in the expectation of \eqref{HN} only depends
on $l_{N}$ and the position
of the return to
the origin, namely $i_{1},...,i_{l_{N}}$. Hence we can rewrite $H_{N}$ as an expectation
 over $P_{\alpha,h}^{\beta}$ since we know the Radon Nikodym density $dP
 /dP_{\alpha,h}^{\beta}(\{i_{1},...,i_{l_{N}}\})$.
Hence $H_{N}$ becomes
\begin{align*}
H_{N}&=E_{\alpha,h}^{\beta}\left(\exp\left(\ \beta s\ \sum_{j=1}^{l_{N}}
\ \zeta_{i_{j}}\right)\prod_{j=1}^{l_{N}}\ \frac{H_{\alpha,h}^{\beta}}{\alpha^{\tau_{j}}}\ \left(\frac{1+
e^{-2h\left(N-i_{l_{N}}
\right)}}{2}\right) \frac{P\left(\tau\geq N-i_{l_{N}}\right)}{P_{\alpha,h}^{\beta}\left(\tau\geq N-i_{l_{N}}\right)
}\right)
\end{align*}
Now we aim at transforming the excursions law again, so that the chain comes back more often in sites
 where the pinning reward is large. In fact we want the chain to take into account its local environment.
 So we define $P_{\alpha,h}^{\beta,\zeta,\alpha_{1}}$ the law of a non homogenous Markov process
 which depends on the environment. Its excursion laws are defined as follow.
 We set:
$\alpha_{1}<\left(1-P_{\alpha,h}^{\beta}\left(\tau=2\right)\right)/P_{\alpha,h}^{\beta}\left(\tau=2\right)$, such that $\
\mu_{1}=1-\left(\alpha_{1} P_{\alpha,h}^{\beta}\left(\tau=2\right)\right)/
\left(1-P_{\alpha,h}^{\beta}\left(\tau=2\right)\right)>0$ and
\medskip

\begin{align}\label{mes2}
\nonumber &P_{\alpha,h}^{\beta,\zeta,\alpha_{1}}\left(\tau=2\right)=P_{\alpha,h}^{\beta}\left(\tau=2\right)
\left(1+\alpha_{1}
\right)^{\boldsymbol{1}_{\{\zeta_{2}>0\}}}\\
&P_{\alpha,h}^{\beta,\zeta,\alpha_{1}}\left(\tau=2r\right)=P_{\alpha,h}^{\beta}\left(\tau=2r\right)\mu_{1}
^{\boldsymbol{1}_{\{\zeta_{2}>0\}}}\ \,\text{for}\ r\geq2
\end{align}
So, under the law of this process, if the chain comes back to the origin at time i, the law of the following
 excursion is
$P_{\alpha,h}^{\beta,\zeta_{i+.},\alpha_{1}}$. Thus the chain checks wether
the reward at time $i+2$ is positive or negative. If $\zeta_{i+2}\geq 0$ the probability to come back
 to zero at time $i+2$ increases. Else it remains the same.

 With this new process we can write
\begin{multline*}
 H_{N}=E_{\alpha,h}^{\beta,\zeta,\alpha_{1}}\left(\exp\left(\ \beta s\, \sum_{j=1}^{l_{N}}
\zeta_{i_{j}}\right)
\,\prod_{j=1}^{l_{N}}\left(\frac{H_{\alpha,h}^{\beta}}{\alpha^{\tau_{j}}}\right)
\left(\frac{1}{2}+
\frac{e^{-2h\left(N-i_{l_{N}}\right)}}{2}\right)\right.\\\left.
\prod_{j=1}^{l_{N}}
\left(\frac{P_{\alpha,h}^{\beta}\left(\tau_{j}\right)}{P_{\alpha,h}^{\beta,\zeta_{i_{j-1}+.},\alpha_{1}}
\left(\tau_{j}\right)}\right)\
\frac{P\left(\tau\geq N-i_{l_{N}}\right)}{P_{\alpha,h}^{\beta,\zeta_{i_{l_{N}}+.},\alpha_{1}}
\left(\tau\geq N-i_{l_{N}}\right)}\right)
\end{multline*}
\begin{equation*}
H_{N}\geq E_{\alpha,h}^{\beta,\zeta,\alpha_{1}}\left(\exp\left(\beta s\sum_{j=1}^{l_{N}}
\zeta_{i_{j}}\right)
\left(H_{\alpha,h}^{\beta}\right)^{l_{N}}\right.\\\left.
\frac{1}{2}\ \prod_{j=1}^{l_{N}}
\left(\frac{P_{\alpha,h}^{\beta}\left(\tau_{j}\right)}{P_{\alpha,h}^{\beta,\zeta_{i_{j-1}+.},\alpha_{1}}
\left(\tau_{j}\right)}\right)
P\left(\tau\geq N-i_{l_{N}}\right)\right)
\medskip
\end{equation*}
Now we apply the Jensen formula and
\begin{align}\label{eq:min}
\mathbb{E}\left(\frac{1}{N}\log H_{N}\right)&\geq \frac{\beta s}{N}\ \mathbb{E}E_{\alpha,h}
^{\beta,\zeta,\alpha_{1}}
\left(\sum_{j=1}^{l_{N}}
\zeta_{i_{j}}\right)
+\log\left(H_{\alpha,h}^{\beta}\right) \mathbb{E}E_{\alpha,h}^{\beta,\zeta,\alpha_{1}}\left(\frac{l_{N}}{N}\right)
+\frac{1}{N}\log\left(\frac{1}{2}\right)\\ \nonumber
 &+\frac{1}{N}\mathbb{E}E_{\alpha,h}^{\beta,\zeta,\alpha_{1}}\left(\sum_{j=1}^{l_{N}}
\log\left(\frac{P_{\alpha,h}^{\beta}\left(\tau_{j}\right)}{P_{\alpha,h}^{\beta,\zeta_{i_{j-1}+.},\alpha_{1}}
\left(\tau_{j}\right)}\right)\right)+\frac{1}{N}\log\left(
P\left(\tau\geq N\right)\right)
\end{align}
At this point, we can divide in two parts the lower bound of \eqref{eq:min}. The first one (called $E_{1}(N)$)
is a positive energetic
term corresponding to the additional reward the chain can expect on coming back often
in "high reward" sites. Namely
\begin{equation*}
E_{1}(N)=\frac{\beta s}{N}\ \mathbb{E}E_{\alpha,h}^{\beta,\zeta,\alpha_{1}}
\left(\sum_{j=1}^{l_{N}}
\zeta_{i_{j}}\right)
\end{equation*}
the second one ($E_{2}(N)$)is a negative entropic term, because the measures transformations we did
have an entropic cost, namely
\begin{align*}
E_{2}(N)&=\log\left(H_{\alpha,h}^{\beta}\right) \mathbb{E}E_{\alpha,h}^{\beta,\zeta,\alpha_{1}}
\left(\frac{l_{N}}{N}\right)
+\frac{1}{N}\log\left(\frac{1}{2}\right)\\
 &+\frac{1}{N}\mathbb{E}E_{\alpha,h}^{\beta,\zeta,\alpha_{1}}\left(\sum_{j=1}^{l_{N}}
\log\left(\frac{P_{\alpha,h}^{\beta}\left(\tau_{j}\right)}{P_{\alpha,h}^{\beta,\zeta_{i_{j-1}+.},\alpha_{1}}
\left(\tau_{j}\right)}\right)\right)+\frac{1}{N}\log\left(
P\left(\tau\geq N\right)\right)
\end{align*}

\subsubsection{STEP2: energy term computation}

First remark that
\begin{equation}\label{zeta}
\sum_{j=1}^{l_{N}}\zeta_{i_{j}}=\sum_{i=0}^{N-2} \zeta_{i+2}\, \boldsymbol{1}_{\{S_{i}=0\}}
\,\boldsymbol{1}_{\{S_{i+2}=0\}}
+\sum_{k=3}^{N}\sum_{s=0}^{N-k}\zeta_{s+k}\, \boldsymbol{1}_{\{S_{s}=0\}}\, \boldsymbol{1}_{\{S_{i}\neq0\
\forall i \in
\{s+1,..,s+k-1\}\ \text{and} S_{s+k}=0\}}
\end{equation}
\bigskip
 So we put $A=\sum_{i=0}^{N-2} \zeta_{i+2}\, \boldsymbol{1}_{\{S_{i}=0\}}
\,\boldsymbol{1}_{\{S_{i+2}=0\}}$\\
\bigskip and\ \ \  $B=\sum_{k=3}^{N}\sum_{s=0}^{N-k}\zeta_{s+k}\,
 \boldsymbol{1}_{\{S_{s}=0\}}\, \boldsymbol{1}_{\{S_{i}\neq0\
\forall i \in
\{s+1,..,s+k-1\}\ \text{and} S_{s+k}=0\}}$\\
Hence we can compute separately the contributions of $A$ and $B$
\begin{equation*}
\mathbb{E}E_{\alpha,h}^{\beta,\zeta,\alpha_{1}}\left(B\right)=\sum_{k=3}^{N}\sum_{s=0}^{N-k}
\mathbb{E}E_{\alpha,h}^{\beta,\zeta,\alpha_{1}}
\left(\zeta_{s+k}\,
 \boldsymbol{1}_{\{S_{s}=0\}}\, \boldsymbol{1}_{\{S_{i}\neq0\
\forall i \in
\{s+1,..,s+k-1\}\ \text{and} S_{s+k}=0\}}\right)\\
\end{equation*}
By Markov property
\begin{align*}
\mathbb{E}E_{\alpha,h}^{\beta,\zeta,\alpha_{1}}\left(B\right)=
\sum_{k=3}^{N}\sum_{s=0}^{N-k}&
\mathbb{E}\left(\boldsymbol{1}_{\{\zeta_{s+2}>0\}}E_{\alpha,h}^{\beta,\zeta,\alpha_{1}}\left(
 \boldsymbol{1}_{\{S_{s}=0\}}\right) P_{\alpha,h}^{\beta}\left(k\right)\mu_{1}\,\zeta_{s+k}\right)\\
 +&\mathbb{E}\left(\boldsymbol{1}_{\{\zeta_{s+2}\leq0\}}E_{\alpha,h}^{\beta,\zeta,\alpha_{1}}\left(
 \boldsymbol{1}_{\{S_{s}=0\}}\right) P_{\alpha,h}^{\beta}\left(k\right)\,\zeta_{s+k}\right)\\
 \end{align*}
 But we notice that $E_{\alpha,h}^{\beta,\zeta,\alpha_{1}}\left(
 \boldsymbol{1}_{\{S_{s}=0\}}\right)$ only depends on $\{\zeta_{1},\zeta_{2},...,\zeta_{s}\}$,
  hence by independence of the $\{\zeta_{i}\}_{i\geq 1}$ and since they are centered and $k\geq3$ we have:
\
$\mathbb{E}E_{\alpha,h}^{\beta,\zeta,\alpha_{1}}\left(B\right)=0$.

Now let's consider the contribution of part A in \eqref{zeta}
\begin{align*}
\mathbb{E}E_{\alpha,h}^{\beta,\zeta,\alpha_{1}}\left(A\right)=&\sum_{i=0}^{N-2}
\mathbb{E}\left(E_{\alpha,h}^{\beta,\zeta,\alpha_{1}}
\left(
 \boldsymbol{1}_{\{S_{i}=0\}}\right)\, P_{\alpha,h}^{\beta}\left(2\right) \left(1+\alpha_{1}\right)
 \zeta_{i+2} \boldsymbol{1}_{\{\zeta_{i+2}>0\}}\right)\\
 &+\sum_{i=0}^{N-2}
\mathbb{E}\left(E_{\alpha,h}^{\beta,\zeta,\alpha_{1}}
\left(
 \boldsymbol{1}_{\{S_{i}=0\}}\right)\, P_{\alpha,h}^{\beta}\left(2\right)
 \zeta_{i+2} \boldsymbol{1}_{\{\zeta_{i+2}\leq0\}}\right)\\
 =&\alpha_{1}P_{\alpha,h}^{\beta}\left(2\right)\mathbb{E}
 \left(\zeta_{1}\,\boldsymbol{1}_{\{\zeta_{1}\}>0}\right)
 \mathbb{E}E_{\alpha,h}^{\beta,\zeta,\alpha_{1}}\left(\sharp\{i\in\{0,..,N-2\} : S_{i}=0\}\right)
\end{align*}
\noindent So the contribution of this energy term is
\begin{equation}\label{E1}
E_{1}(N)=\beta s \alpha_{1}P_{\alpha,h}^{\beta}\left(2\right)
 \frac{\mathbb{E}E_{\alpha,h}^{\beta,\zeta,\alpha_{1}}\left(\sharp\{i\in\{0,..,N-2\} :
 S_{i}=0\}\right)}{N}\geq \beta s \alpha_{1}P_{\alpha,h}^{\beta}\left(2\right)
 \frac{\mathbb{E}E_{\alpha,h}^{\beta,\zeta,\alpha_{1}}\left(l_{N}\right)}{N}
\end{equation}
\subsubsection{STEP3: computation of entropic term}
First notice that the terms $1/N\log\left(
P\left(\tau\geq N\right)\right)$ and $1/N \log(1/2)$ go to $0$ as $N$
goes to $\infty$ independently of all the other parameters. So we put $R_{N}=1/N \log\left(
P\left(\tau\geq N\right)\right)+ 1/N \log\left(1/2\right)$ and we can write
\begin{equation*}
E_{2}(N)=\frac{S_{N}}{N}+\log\left(H_{\alpha,h}^{\beta}\right) \mathbb{E}E_{\alpha,h}^{\beta,\zeta,\alpha_{1}}
\left(\frac{l_{N}}{N}\right)+R_{N}
\end{equation*}
where we have put
\begin{equation}
S_{N}=\mathbb{E}E_{\alpha,h}^{\beta,\zeta,\alpha_{1}}\left(\sum_{j=1}^{l_{N}}
\log\left(\frac{P_{\alpha,h}^{\beta}\left(\tau_{j}\right)}{P_{\alpha,h}^{\beta,\zeta_{i_{j-1}+.},\alpha_{1}}
\left(\tau_{j}\right)}\right)\right)
\end{equation}
The definitions \eqref{mes} and \eqref{mes2} of $ P_{\alpha,h}^{\beta,\zeta_{i_{j-1}+.},\alpha_{1}}$
 and $P_{\alpha,h}^{\beta}$
give us immediately
\begin{align*}
S_{N}=&-\mathbb{E}E_{\alpha,h}^{\beta,\zeta,\alpha_{1}}\left(\sum_{j=1}^{l_{N}}
 \boldsymbol{1}_{\{\zeta_{i_{j-1}+2}>0\}}\left( \boldsymbol{1}_{\{\tau_{j}=2\}} \log\left(1+\alpha_{1}\right)
+\boldsymbol{1}_{\{\tau_{j}>2\}} \log\left(\mu_{1}\right)\right)\right)\\
=&-\sum_{i=0}^{N-2}\mathbb{E}\left(E_{\alpha,h}^{\beta,\zeta,\alpha_{1}}\left(\boldsymbol{1}_{\{S_{i}=0\}}
\boldsymbol{1}_{\{S_{i+2}=0\}}\right)
 \boldsymbol{1}_{\{\zeta_{i+2}>0\}} \log\left(1+\alpha_{1}\right)\right)\\
&-\sum_{k=3}^{N}\sum_{s=0}^{N-k}\mathbb{E}\left(E_{\alpha,h}^{\beta,\zeta,\alpha_{1}}\left(
\boldsymbol{1}_{\{S_{s}=0\}}\boldsymbol{1}_{\{S_{s+k}=0\}}\boldsymbol{1}_{\{S_{i}\neq0\
\forall i \in
\{s+1,..,s+k-1\}\}}\right)\boldsymbol{1}_{\{\zeta_{s+2}>0\}}
\log\left(\mu_{1}\right)\right)\\
\end{align*}
And once again, by Markov property we have
\begin{equation*}
 \boldsymbol{1}_{\{\zeta_{i+2}>0\}}\ E_{\alpha,h}^{\beta,\zeta,\alpha_{1}}\left(\boldsymbol{1}_{\{S_{i}=0\}}
\boldsymbol{1}_{\{S_{i+2}=0\}}\right)=\boldsymbol{1}_{\{\zeta_{i+2}>0\}}\ E_{\alpha,h}^{\beta,\zeta,\alpha_{1}}\left(
\boldsymbol{1}_{\{S_{i}=0\}}\right) \left(1+\alpha_{1}\right)P_{\alpha,h}^{\beta}\left(2\right)
\end{equation*}
We notice that $E_{\alpha,h}^{\beta,\zeta,\alpha_{1}}\left(\boldsymbol{1}_{\{S_{i}=0\}}\right)$
is independant of $\zeta_{i+2}$ and $\mathbb{P}\left(\zeta_{i+2}>0\right)=1/2$ hence
\begin{align*}
S_{N}=&-\frac{P_{\alpha,h}^{\beta}\left(2\right)}{2}
\left(1+\alpha_{1}\right) \log\left(1+\alpha_{1}\right)
\mathbb{E}E_{\alpha,h}^{\beta,\zeta,\alpha_{1}}\left(l_{N-2}\right)\\
&-\sum_{k=3}^{N}\ \frac{\mu_{1} \log\left(\mu_{1}\right)}{2}\
P_{\alpha,h}^{\beta}\left(k\right)\
\mathbb{E}E_{\alpha,h}^{\beta,\zeta,\alpha_{1}}\left(
l_{N-k}\right)
\end{align*}
Finally the entropic contribution is
\begin{align} \label{E2}
\nonumber E_{2}(N)&=\log\left(H_{\alpha,h}^{\beta}\right)
\mathbb{E}E_{\alpha,h}^{\beta,\zeta,\alpha_{1}}
\left(\frac{l_{N}}{N}\right)-\frac{1}{2}
P_{\alpha,h}^{\beta}\left(2\right) \left(1+\alpha_{1}\right)
\log\left(1+\alpha_{1}\right)
\mathbb{E}E_{\alpha,h}^{\beta,\zeta,\alpha_{1}}\left(\frac{l_{N-2}}{N}\right)\\
&-\sum_{k=3}^{N}\ \frac{\mu_{1} \log\left(\mu_{1}\right)}{2}\
P_{\alpha,h}^{\beta}\left(k\right)\
\mathbb{E}E_{\alpha,h}^{\beta,\zeta,\alpha_{1}}\left(
\frac{l_{N-k}}{N}\right)+R_{N}
\end{align}

So \eqref{E1} and \eqref{E2} give us a precise lower bound of formula
\eqref{eq:min} of the form
\begin{equation}\label{min2}
 \mathbb{E}\left(\frac{1}{N}\log \left(H_{N}\right)\right)\geq
E_{1}(N)+E_{2}(N)
\end{equation}

\subsubsection{STEP4: estimation of $H_{\alpha,h}^{\beta}$ and choice of $\alpha$ and $\alpha_{1}$}
Now we want to evaluate $H_{\alpha,h}^{\beta}$ with its expression of \eqref{mes}
\begin{equation*}
H_{\alpha,h}^{\beta}=e^{\beta}\left(1-\frac{\sqrt{1-\alpha^{2}}+\sqrt{1-e^{-4h}\alpha^{2}}}{2}\right)\\
\end{equation*}
In order to compare $\log\left(H_{\alpha,h}^{\beta}\right)$ with the other terms of \eqref{min2}, we
put $\alpha^{2}=1-c \alpha_{1}^{2}$, with $c>0$ and $\sqrt{c} \alpha_{1} \leq 1$. That way we obtain
\begin{align*}
H_{\alpha,h}^{\beta}=&e^{\beta}\left(1-\frac{\sqrt{1-e^{-4h}}}{2}+
\frac{\sqrt{1-e^{-4h}}-\sqrt{1-e^{-4h}\left(1-c\alpha_{1}^{2}\right)}-\sqrt{c}\alpha_{1}}{2}\right)\\
H_{\alpha,h}^{\beta}=&e^{\beta}\left(1-\frac{\sqrt{1-e^{-4h}}}{2}\right)
\left(1+\frac{\sqrt{1-e^{-4h}}\left(1-\sqrt{1+\frac{c e^{-4h}\alpha_{1}^
{2}}{1-e^{-4h}}}\right)-\sqrt{c}\alpha_{1}}{2-\sqrt{1-e^{-4h}}}\right)
\end{align*}
But $\sqrt{1+x}\leq 1+ x/2$ for $x \in (-1,+ \infty)$ and
$2-\sqrt{1-e^{-4h}}\geq 1$ hence:
\begin{equation*}
\log\left(H_{\alpha,h}^{\beta}\right)\geq \log \left(e^{\beta}\left(1-\frac{\sqrt{1-e^{-4h}}}{2}\right)\right)+
\log \left(1-\sqrt{c}\alpha_{1}-\frac{c \alpha_{1}^{2}e^{-4h}}{2\sqrt{1-e^{-4h}}}\right)
\end{equation*}
But as $\sqrt{c}\alpha_{1}\leq 1$ we can bound by above the term
\begin{equation}\label{inter}
\sqrt{c}\alpha_{1}+\frac{c
\alpha_{1}^{2}e^{-4h}}{2\sqrt{1-e^{-4h}}}=\sqrt{c}
\alpha_{1}\left(1+\frac{\sqrt{c}\alpha_{1}
e^{-4h}}{2\sqrt{1-e^{-4h}}}\right)
\leq \sqrt{c}\alpha_{1}\left(1+\frac{
1}{2\sqrt{1-e^{-4h}}}\right)\\
\end{equation}
To continue our computation we need to choose precise
values for $\alpha_{1}$ and $c$. That is why recalling that
$\left(\ \alpha^{2}=1-c\alpha_{1}^{2}\right)$ we put
\begin{equation}\label{alpha1}
 \alpha_{1}=\beta s/\left(5\times 2^{8}\right)\ \ \ \
 \sqrt{c}=\beta s/\left(3\times 2^{4}\left(1+\frac{1}{2\sqrt{1-e^{-4h}}}\right)\right)
\end{equation}
Notice that $\log(1-x)\geq -3x/2$ if $x \in [0,1/3]$, and since
$\beta s\leq \log(2)$ the rhs of \eqref{inter} verify
$\sqrt{c}\alpha_{1}\left(1+1/\left(2\sqrt{1-e^{-4h}}\right)\right)\leq
\beta^{2}s^{2}/\left(15 \times 2^{12}\right)\leq
 \frac{1}{3}$ hence $\log\left(H_{\alpha,h}^{\beta}\right)$
becomes

 \begin{align*}
\log\left(H_{\alpha,h}^{\beta}\right)\geq &\log \left(e^{\beta}\left(1-\frac{\sqrt{1-e^{-4h}}}{2}\right)\right)
-\frac{3}{2}\sqrt{c}\alpha_{1}
\left(1+\frac{1}{2\sqrt{1-e^{-4h}
}}\right)\\
\geq &\log \left(e^{\beta}\left(1-\frac{\sqrt{1-e^{-4h}}}{2}\right)\right)
-\frac{\beta^{2}s^{2}}{5\times 2^{13}}
\end{align*}
Hence as $\log(1+\alpha_{1})\leq \alpha_{1}$ we can rewrite equation \eqref{eq:min}
\begin{align}\label{min3}
\nonumber \mathbb{E}\left(\frac{1}{N}\log \left(H_{N}\right)\right)\geq
&\Bigg[\beta s \alpha_{1} P_{\alpha,h}^{\beta}\left(2\right)
-\frac{1}{2}P_{\alpha,h}^{\beta}\left(2\right)\left(1+\alpha_{1}\right)\alpha_{1}\\
\nonumber &+\log\left(e^{\beta}\left(1-\frac{\sqrt{1-e^{-4h}}}{2}\right)\right)
-\frac{\beta^{2}s^{2}}{5\times 2^{13}}
\Bigg]
 \mathbb{E}\left( E_{\alpha,h}^{\beta,\zeta,\alpha_{1}}\left(\frac{l_{N}}{N}\right)\right)\\
&-\sum_{k=3}^{N}P_{\alpha,h}^{\beta}\left(k\right)\frac{\mu_{1}\log\left(\mu_{1}\right)}{2}
\mathbb{E}\left( E_{\alpha,h}^{\beta,\zeta,\alpha_{1}}\left(\frac{l_{N-k}}{N}\right)\right)+R_{N}
\end{align}
\subsubsection{STEP5: intermediate computation}
To conclude this computation we need some inequalities on $P_{\alpha,h}^{\beta}$ and
$H_{\alpha,h}^{\beta}$ .
As $\beta s \leq \log(2)$ equations \eqref{alpha1} show that $\alpha_{1} \sqrt{c}\in [0,1/4]$, hence
$\alpha^{2}=1-c\alpha_{1}^{2}\geq 1-1/2^{4}\geq 3/4$. 
So we can bound from above and below the quantity $H_{\alpha,h}^{\beta}$ (introduced in \eqref{H})
\begin{equation*}
e^{\beta}\geq H_{\alpha,h}^{\beta}
\geq
e^{\beta}\left(1-\frac{\sqrt{c}\alpha_{1}}{2}-\frac{1}{2}\right)\geq \frac{3e^{\beta}}{8}
\end{equation*}
At this point we need to bound from above and below the quantity $P_{\alpha,h}^{\beta}\left(2\right)$,
 which has been defined in \eqref{mes}. With the previous inequalities we have $e^{\beta}/H_{\alpha,h}^{\beta}\geq 1$
 and $\sqrt{1-\alpha^{2}}\leq1/4$ so
\begin{equation}\label{W}
P_{\alpha,h}^{\beta}\left(2\right)=1-\sum_{i=2}^{\infty}P_{\alpha,h}^{\beta}\left(2i\right)\leq
1-\sum_{i=2}^{\infty}\frac{1}{2} \alpha^{2i}
 P(\tau=2i)=1-\frac{1}{2}\left(1-\sqrt{1-\alpha^{2}}-\frac{\alpha^{2}}{2}\right)\leq \frac{7}{8}
\end{equation}
and
\begin{equation}\label{M3}
\frac{1}{8}=\frac{1}{4}\times \frac{e^{\beta}}{2e^{\beta}}\leq
P_{\alpha,h}^{\beta}\left(2\right)
\end{equation}
\noindent
And to finish with these preliminary inequalities, we notice with \eqref{W} and \eqref{M3} that
\begin{equation}\label{M4}
\frac{1}{8}\leq 1-P_{\alpha,h}^{\beta}\left(2\right)\ \ \ \text{and}\ \ \
\frac{1}{7}\leq \frac{P_{\alpha,h}^{\beta}\left(2\right)}{1-P_{\alpha,h}^{\beta}\left(2\right)}\leq
7
\end{equation}
Hence the condition $\alpha_{1}<P_{\alpha,h}^{\beta}\left(\tau=2\right)/
\left(1-P_{\alpha,h}^{\beta}\left(\tau=2\right)\right)$ is obviously verified.
\subsubsection{STEP 6: conclusion}
In the equation \eqref{min3} we still have to evaluate the term
\begin{equation*}
\sum_{k=3}^{N}P_{\alpha,h}^{\beta}\left(k\right)\mathbb{E}\left( E_{\alpha,h}^{\beta,\zeta,\alpha_{1}}
\left(\frac{l_{N-k}}{N}\right)\right)
\end{equation*}
So if $N\geq N_{0}$
\begin{align*}
\sum_{k=3}^{N}P_{\alpha,h}^{\beta}\left(k\right)\mathbb{E}\left( E_{\alpha,h}^{\beta,\zeta,\alpha_{1}}
\left(\frac{l_{N-k}}{N}\right)\right)
\geq &P_{\alpha,h}^{\beta}{\left(\{3,..,N_{0}\}\right)}\mathbb{E} E_{\alpha,h}^{\beta,\zeta,\alpha_{1}}
\left(\frac{l_{N-N_{0}}}{N}\right)\\
\geq &\left(1-P_{\alpha,h}^{\beta}\left(2\right)\right)\mathbb{E} E_{\alpha,h}^{\beta,\zeta,\alpha_{1}}
\left(\frac{l_{N}}{N}\right)-\frac{N_{0}}{N}\\
&-P_{\alpha,h}^{\beta}\left(\{N_{0}+1,..
,\infty\}\right)\mathbb{E}E_{\alpha,h}^{\beta,\zeta,\alpha_{1}}\left(\frac{l_{N}}{N}\right)
\end{align*}
\noindent
Hence equation \eqref{min3} becomes
\begin{align}\label{min4}
\nonumber \mathbb{E}\left(\frac{1}{N}\log \left(H_{N}\right)\right)\geq
&\Bigg[\beta s \alpha_{1} P_{\alpha,h}^{\beta}\left(2\right)
-\frac{1}{2}P_{\alpha,h}^{\beta}\left(2\right)\left(1+\alpha_{1}\right)\alpha_{1}
-\frac{\beta^{2}s^{2}}{5\times 2^{13}}\\ \nonumber
&+\log\left(e^{\beta}\left(1-\frac{\sqrt{1-e^{-4h}}}{2}\right)\right)-\left(1-P_{\alpha,h}^{\beta}\left(2\right)
\right)
\frac{\mu_{1}\log\left(\mu_{1}\right)}{2}\\ \nonumber
&+P_{\alpha,h}^{\beta}\left(\{N_{0}+1,..
,\infty\}\right)\frac{\mu_{1}\log\left(\mu_{1}\right)}{2}\Bigg]
\mathbb{E}\left( E_{\alpha,h}^{\beta,\zeta,\alpha_{1}}\left(\frac{l_{N}}{N}\right)\right)\\
&+\frac{N_{0}}{N}\frac{\mu_{1}\log\left(\mu_{1}\right)}{2}+R_{N}
\end{align}
We can now bound from below with \eqref{alpha1} and \eqref{M3}
\begin{equation*}
\beta s \alpha_{1}P_{\alpha,h}^{\beta}\left(2\right)\geq \frac{\beta s}{2^{3}}\frac{\beta s}{5\times
2^{8}}
= \frac{\beta^{2} s^{2}}{5\times 2^{11}}
\end{equation*}
Moreover
$\mu_{1}=1-\frac{\alpha{1}P_{\alpha,h}^{\beta}\left(2\right)}{1-P_{\alpha,h}^{\beta}\left(2\right)}$\
\  and  $-\log(1-x)\geq x$\  for\  $x\in [0,1)$
 so we have
\begin{equation*}
-\frac{1-P_{\alpha,h}^{\beta}\left(2\right)}{2}\mu_{1}\log\left(\mu_{1}\right)\geq
\frac{\alpha_{1}P_{\alpha,h}^{\beta}\left(2\right)}{2}-\frac{\alpha_{1}^{2}P_{\alpha,h}^{\beta}
\left(2\right)^{2}}{2\left(
1-P_{\alpha,h}^{\beta}\left(2\right)\right)}
\end{equation*}
We noticed before in \eqref{M3} and \eqref{M4} that
$P_{\alpha,h}^{\beta}\left(2\right)\leq 7/8$ and
$P_{\alpha,h}^{\beta}\left(2\right)
/\left(2\left(1-P_{\alpha,h}^{\beta}\left(2\right)\right)\right)\leq
\frac{7}{2}$,\  hence
\begin{equation*}
-\frac{1-P_{\alpha,h}^{\beta}\left(2\right)}{2}\mu_{1}\log\left(\mu_{1}\right)\geq
\frac{\alpha_{1}P_{\alpha,h}^{\beta}\left(2\right)}{2}-\frac{7^{2}\alpha_{1}^{2}}{2^{4}}
\geq\frac{\alpha_{1}P_{\alpha,h}^{\beta}\left(2\right)}{2}-4\alpha_{1}^{2}
\end{equation*}
That way the inequality \eqref{min4} must now be written

\begin{align}\label{min5}
\nonumber \mathbb{E}\left(\frac{1}{N}\log
\left(H_{N}\right)\right)\geq &\Bigg[\frac{\beta^{2}
s^{2}}{5\times 2^{12}}
-\frac{1}{2}P_{\alpha,h}^{\beta}\left(2\right)\left(1+\alpha_{1}\right)\alpha_{1}+\frac{\alpha_{1}
P_{\alpha,h}^{\beta}\left(2\right)}{2}-4
\alpha_{1}^{2}\\ \nonumber
&+\log\left(e^{\beta}\left(1-\frac{\sqrt{1-e^{-4h}}}{2}\right)\right)\\
\nonumber &+P_{\alpha,h}^{\beta}\left(\{N_{0}+1,..
,\infty\}\right)\frac{\mu_{1}\log\left(\mu_{1}\right)}{2}\Bigg]
\mathbb{E}\left(
E_{\alpha,h}^{\beta,\zeta,\alpha_{1}}\left(\frac{l_{N}}{N}\right)\right)\\
&+\frac{N_{0}}{N}\mu_{1}\log\left(\mu_{1}\right)+R_{N}
\end{align}
By \eqref{M4} and \eqref{M3} we know that $P_{\alpha,h}^{\beta}(2)
\leq 7/8$ and
$P_{\alpha,h}^{\beta}\left(2\right)/\left(1-P_{\alpha,h}^{\beta}\left(2\right)\right)\leq
7$. Hence we have the inequalities
\begin{align}\label{M5}
-\frac{1}{2}P_{\alpha,h}^{\beta}\left(2\right)&\left(1+\alpha_{1}\right)\alpha_{1}+\frac{\alpha_{1}
P_{\alpha,h}^{\beta}\left(2\right)}{2}-4
\alpha_{1}^{2}\geq -5 \alpha_{1}^{2} \geq -\frac{\beta^{2} s^{2}}{5\times 2^{16}}\\
\text{and}\ \ \ \ \  &\frac{\alpha_{1}
P_{\alpha,h}^{\beta}\left(2\right)}{1-P_{\alpha,h}^{\beta}\left(2\right)}\leq
7 \alpha_{1}= \frac{7 \beta s}{5\times 2^{8}} <\frac{1}{3}
\end{align}
Now since $\mu_{1}\leq 1$ and $\log\left(1-x\right)\geq -3x/2$ for
$x\in [0,1/3]$ equation $2.20$ allows us to bound by below
\begin{equation*}
\mu_{1}\log\left(\mu_{1}\right)\geq -\frac{3}{2}\ \frac{P_{\alpha,h}^{\beta}\left(2\right)}
{1-P_{\alpha,h}^{\beta}\left(2\right)}\ \alpha_{1} \geq -\frac{21 \beta
s}{5\times 2^{9}}\geq -1
\end{equation*}
So equation \eqref{min5} becomes
\begin{align} \label{min6}
\nonumber\mathbb{E}\left(\frac{1}{N}\log
\left(H_{N}\right)\right)\geq
&\Bigg[\frac{\beta^{2} s^{2}}{5\times 2^{13}}
+\log\left(e^{\beta}\left(1-\frac{\sqrt{1-e^{-4h}}}{2}\right)\right)\\
&-P_{\alpha,h}^{\beta}\left(\{N_{0}+1,..
,\infty\}\right)
\Bigg] \mathbb{E}\left(
E_{\alpha,h}^{\beta,\zeta,\alpha_{1}}\left(\frac{l_{N}}{N}\right)\right)
-\frac{N_{0}}{N}
+R_{N}
\end{align}
But as proved in appendix A.$1$), $P_{\alpha,h}^{\beta}\left(\{N_{0}+1,..
,\infty\}\right)$ goes to zero as $N_{0}$ goes to
$\infty$ independently of $h \geq 0$,
hence for $N_{0}$ large enough and for all $h>0$
\begin{equation*}
P_{\alpha,h}^{\beta}\left(\{N_{0}+1,..
,\infty\}\right)\leq \frac{\beta^{2} s^{2}}{5\times 2^{14}}
\end{equation*}
So if we put $q\left(s\right)=\frac{\beta^{2}
s^{2}}{5\times2^{14}}$ the equation \eqref{min6} gives us for all
$N\geq N_{0}$ and $h>0$
\begin{equation}\label{min8}
\mathbb{E}\left(\frac{1}{N}\log \left(H_{N}\right)\right)\geq
\Bigg[q\left(s\right)
+\log\left(e^{\beta}\left(1-\frac{\sqrt{1-e^{-4h}}}{2}\right)\right)\Bigg]
\mathbb{E}\left(
E_{\alpha,h}^{\beta,\zeta,\alpha_{1}}\left(\frac{l_{N}}{N}\right)\right)+R_{N}^{N_{0}}
\end{equation}
with $R_{N}^{N_{0}}=R_{N}-N_{0}/N$.

 As proved in appendix A.$2$) for every $N\geq 1$
$\mathbb{E}\left(
E_{\alpha,h}^{\beta,\zeta,\alpha_{1}}\left(l_{N}/N\right)\right)\geq
\mathbb{E}\left( E_{\alpha,h}^{\beta}\left(l_{N}/N\right)\right)$.
So if we note $h_{0}\left(\beta\right)$ the quantity verifying
$\log\left(e^{\beta}\left(1-\sqrt{1-e^{-4h_{o}\left(\beta\right)}}/2\right)\right)=-q\left(s\right)$
we have for every $h<h_{0}\left(\beta\right)$ and $N\geq N_{0}$ that
\begin{equation*}
\mathbb{E}\left(\frac{1}{N}\log \left(H_{N}\right)\right)\geq
\Bigg[q\left(s\right)
+\log\left(e^{\beta}\left(1-\frac{\sqrt{1-e^{-4h}}}{2}\right)\right)\Bigg]
\ \mathbb{E}\left(
E_{\alpha,h}^{\beta}\left(\frac{l_{N}}{N}\right)\right)+R_{N}^{N_{0}}
\end{equation*}
and consequently
\begin{equation*}
\Phi^{s}\left(\beta,h\right)\geq
\Bigg[q\left(s\right)
+\log\left(e^{\beta}\left(1-\frac{\sqrt{1-e^{-4h}}}{2}\right)\right)\Bigg]
\ \liminf_{N \to \infty}\mathbb{E}\left( E_{\alpha,h}^{\beta}\left(\frac{l_{N}}{N}\right)\right)
\end{equation*}
Notice also that
 $\liminf_{N \to \infty}\mathbb{E}\left( E_{\alpha,h}^{\beta}\left(\frac{l_{N}}{N}\right)\right)>0$ (because
 $\alpha\in(0,1)$).
Hence for every $\beta$ in  $[0,\log\left(2\right)-q_{s})$, \ $h_{0}(\beta)$ is a lower bound of $h_{c}
\left(\beta\right)$
\begin{equation*}
h_{c}\left(\beta\right)\geq h_{0}\left(\beta\right)=- \frac{1}{4}\log\left(1-4\left(1-e^{-\beta-q\left(s\right)}\right)^{2}\right)
\end{equation*}
\qed

\noindent
\subsection{\textit{Proof of Corollary \ref{cor}}}
As showed just before in \eqref{min8} we have a rank $N_{0}\in \mathbb{N}-\{0\}$ such that for all $h>0$
and $N\geq N_{0}$
\begin{multline*}
\mathbb{E}\left(\frac{1}{N} \log E\left(\exp\left(\beta \sum_{i=1}^{N}\boldsymbol{1}_{\{S_{i}=0\}}\left(s
\zeta_{i}+1\right)
-2h \sum_{i=1}^{N} \Delta_{i}\right)\right)\right)\geq\\
\Bigg[\frac{\beta^{2} s^{2}}{5\times 2^{14}}
+\log\left(e^{\beta}\left(1-\frac{\sqrt{1-e^{-4h}}}{2}\right)\right)\Bigg]
\mathbb{E}\left( E_{\alpha,h}^{\beta,\zeta,\alpha_{1}}\left(\frac{l_{N}}{N}\right)\right)+R_{N}^{N_{0}}
\end{multline*}
but in appendix A.$2$) we prove the following inequalities
\begin{equation}\label{coupl}
\mathbb{E}\left( E_{\alpha,h}^{\beta,\zeta,\alpha_{1}}\left(\frac{l_{N}}{N}\right)\right)\geq
\mathbb{E}\left( E_{\alpha,h}^{\beta}\left(\frac{l_{N}}{N}\right)\right)\geq
\mathbb{E}\left( E_{\alpha,\infty}^{0}\left(\frac{l_{N}}{N}\right)\right)>0
\end{equation}
and for fixed $\beta,s,N$ let $h$ go to $\infty$
\begin{multline*}
\mathbb{E}\left(\frac{1}{N} \log E\left(\exp\left(\beta \sum_{i=1}^{N}\boldsymbol{1}_{\{S_{i}=0\}}\left(s
\zeta_{i}+1\right)
\right)\boldsymbol{1}_{\{S_{i}\geq 0, \forall i\in\{1,..,N\}\}}\right)\right)\geq\\
\Bigg[\frac{\beta^{2} s^{2}}{5\times 2^{14}}
+\log\left(e^{\beta}\frac{1}{2}\right)\Bigg]
\mathbb{E}\left( E_{\alpha,\infty}^{0}\left(\frac{l_{N}}{N}\right)\right)+R_{N}^{N_{0}}
\end{multline*}
Now recall that $P\left(\{S_{i}\geq 0, \forall i\in\{1,..,N\}\}\right)\sim c/N^{1/2}$,
 the lower bound becomes
\begin{multline*}
\mathbb{E}\left(\frac{1}{N} \log E\left(\exp\left(\beta \sum_{i=1}^{N}\boldsymbol{1}_{\{S_{i}=0\}}\left(s
\zeta_{i}+1\right)
\right)\bigg| \{S_{i}\geq 0, \forall i\in\{1,..,N\}\}\right)\right)\geq\\
\Bigg[\frac{\beta^{2} s^{2}}{5\times 2^{14}}
+\log\left(e^{\beta}\frac{1}{2}\right)\Bigg]
\mathbb{E}\left( E_{\alpha,\infty}^{0}\left(\frac{l_{N}}{N}\right)\right)+K_{N}^{N_{0}}
\end{multline*}
With $K_{N}^{N_{0}}=R_{N}^{N_{0}}-1/N \log\left(P\left(\{S_{i}\geq
0, \forall i\in\{1,..,N\}\}\right)\right)$, so that it goes to $0$
as $N$ goes to $\infty$ independently of all the other parameters.
Now by \cite{Yosh} we can apply the fact that for an odd number of
steps the RW conditioned to stay positive becomes the reflected RW
if it is pinned by $\log{2}$, that is to say
\begin{equation*}
\frac{P_{refl. RW}}{P_{RW cond. to be \geq
o}}\left(S\right)=
\frac{e^{\log\left(2\right)\sum_{i=1}^{2N+1}
\boldsymbol{1}_{\{S_{i}=0\}}}\ \ \boldsymbol{1}_{\{S_{i}\geq 0\ \forall i\in\{0,2N+1\}\}}}{V_{2N+1}}
\end{equation*}
With $\frac{1}{N}\log(V_{N})$ goes to $0$ as $N$ goes to $\infty$.
Hence we put $\beta=\log(2)-u$
\begin{multline*}
\mathbb{E}\left(\frac{1}{2N+1} \log E\left(\exp\left(\log(2)
\sum_{i=1}^{2N+1}\boldsymbol{1}_{\{S_{i}=0\}}+
\sum_{i=1}^{2N+1}\boldsymbol{1}_{\{S_{i}=0\}}(-u+\beta s
\zeta_{i})
\right)\bigg| \{S_{i}\geq 0,\forall i \leq 2N+1\}\}\right)\right)\geq\\
\Bigg[\frac{\beta^{2} s^{2}}{5\times 2^{14}} -u\Bigg]
\mathbb{E}\left(
E_{\alpha,\infty}^{0}\left(\frac{l_{2N+1}}{2N+1}\right)\right)+K_{2N+1}^{N_{0}}
\end{multline*}

\begin{multline*}
\mathbb{E}\left(\frac{1}{2N+1} \log E\left(\exp\left(
\sum_{i=1}^{2N+1}\boldsymbol{1}_{\{S_{i}=0\}}(-u+\beta s
\zeta_{i})
\right)\right)\right)\geq\\
\Bigg[\frac{\beta^{2} s^{2}}{5\times 2^{14}} -u\Bigg]
\mathbb{E}\left(
E_{\alpha,\infty}^{0}\left(\frac{l_{2N+1}}{2N+1}\right)\right)+K_{2N+1}^{N_{0}}+\frac{1}{2N+1}\log(V_{2N+1})
\end{multline*}
Now let $N \to \infty$, and recall $\beta=\log(2)-u$
\begin{equation*}
\lim_{N\to\infty}\mathbb{E}\left(\frac{1}{N} \log E\left(\exp\left(
\sum_{i=1}^{N}\boldsymbol{1}_{\{S_{i}=0\}}(-u+\beta s
\zeta_{i})
\right)\right)\right)\geq\Bigg[\frac{\beta^{2} s^{2}}{5\times 2^{14}}
-u\Bigg]
\lim_{N\to\infty} E_{\alpha,\infty}^{0}\left(\frac{l_{N}}{N}\right)
\end{equation*}
so, for $u\leq \log(2)/2$, \  $\left(\beta \geq \log(2)/2\right)$
\begin{equation*}
\lim_{N\to\infty}\mathbb{E}\left(\frac{1}{N} \log E\left(\exp\left(
\sum_{i=1}^{N}\boldsymbol{1}_{\{S_{i}=0\}}(-u+\beta s
\zeta_{i})
\right)\right)\right)\geq\Bigg[\frac{\log(2)^{2} s^{2}}{5\times 2^{16}}
-u\Bigg]
\lim_{N\to\infty} E_{\alpha,\infty}^{0}\left(\frac{l_{N}}{N}\right)
\end{equation*}
By convexity, the free energy $\Phi$, defined by
\begin{equation*}
 \Phi(u,v)=\lim_{N\to\infty}\mathbb{E}\left(\frac{1}{N} \log E\left(\exp\left(
\sum_{i=1}^{N}\boldsymbol{1}_{\{S_{i}=0\}}(-u+v
\zeta_{i})
\right)\right)\right)
\end{equation*}
is not decreasing in $v$ hence
\begin{equation*}
 \Phi(u,\log(2)s)\geq \Bigg[\frac{\log(2)^{2} s^{2}}{5\times 2^{16}}
-u\Bigg]
\lim_{N\to\infty} E_{\alpha,\infty}^{0}\left(\frac{l_{N}}{N}\right)
\end{equation*}
and
for $s \in [0,\log(2)]$
\begin{equation*}
 u_{c}(s)\geq \frac{s^{2}}{5\times 2^{16}}
\end{equation*}
\qed
\appendix

\section{}
\subsection{}
First we have to prove the first point, namely $P_{\alpha,h}^{\beta}(\{N_{0},...,+ \infty\})$ goes to
$0$ as $N_{0}$ goes to infinity independently of $h\geq0$. That way we bound by above the quantity \eqref{mes}
\begin{align*}
P_{\alpha,h}^{\beta}\left(\tau_{1}=2n\right)&=\left(\frac{1+\exp\left(-4hn\right)}{2}\right)\alpha^{2n}
\frac{P\left(\tau=2n\right)}{H_{\alpha,h}^{\beta}}\exp\left(\beta\right)\\
&\leq  \frac{ \alpha^{2n} P\left(\tau=2n\right)}{\sum_{j=1}^{+\infty}\frac{1}{2}\alpha^{2j} P\left(\tau=2j\right)}
\end{align*}
So the rhs of this inequality does not depend on $h$ any more and is the general term of a convergent serie hence
we have the uniform  convergence in $h$.

\subsection{}
Now we want to prove the inequalities of \eqref{coupl}, that is to say
\begin{equation}\label{ineq}
\mathbb{E}\left( E_{\alpha,h}^{\beta,\zeta,\alpha_{1}}\left(\frac{l_{N}}{N}\right)\right)\geq
\mathbb{E}\left( E_{\alpha,h}^{\beta}\left(\frac{l_{N}}{N}\right)\right)\geq
\mathbb{E}\left( E_{\alpha,\infty}^{0}\left(\frac{l_{N}}{N}\right)\right)
\end{equation}
That way we recall a coupling theorem (see \cite{Lig} or \cite{Lin})
\begin{theorem} \label{theo2}
$\mu_{1}$ and $\mu_{2}$ are two probability measures on $2 \mathbb{N}-\{0\}$. If for every bounded and
non decreasing function
$f$ defined on $2 \mathbb{N}-\{0\}$ we have $\mu_{1}(f) \leq \mu_{2}(f)$  we can define on the same
 probability space $\left(\Omega,P\right)$
 two random variables $\left(T_{1},T_{2}\right)$ of law $\left(\mu_{1},\mu_{2}\right)$ such that
$T_{1}\leq T_{2}$ P almost surly.
\end{theorem}
\begin{rem}\label{remark}
We notice that to satisfy the hypothesis of the theorem
it is enough to show that there exists an integer $i_{0}$ such that
$\mu_{1}(2i)\geq \mu_{2}(2i)$ for every $i \in \{1,..,i_{0}\}$ and
$\mu_{1}(2i)\leq \mu_{2}(2i)$ for every $i \geq i_{0}+1$. We can prove it easily on writing
\begin{equation*}
\mu_{2}(f)-\mu_{1}(f)=\sum_{i=1}^{i_{0}}(\mu_{2}(2i)-\mu_{1}(2i)) f(2i)+\sum_{i=i_{0}+1}^{\infty}
(\mu_{2}(2i)-\mu_{1}(2i)) f(2i)
\end{equation*}
But as $f$ is not decreasing $f(2i)\geq f(2i_{0})$ for every $i\geq i_{o}+1$ and $f(2i)\leq f(2i_{0})$
for every $i\leq i_{o}$. Moreover since $\mu_{2}(2i)-\mu_{1}(2i)$ is positive when  $i\geq i_{o}+1$ and
negative else we have the inequality
\begin{align*}
\mu_{2}(f)-\mu_{1}(f)&\geq f(2i_{0})\ \sum_{i=1}^{i_{0}}\mu_{2}(2i)-\mu_{1}(2i)
\ +f(2i_{0})\ \sum_{i=i_{0}+1}^{\infty}
\mu_{2}(2i)-\mu_{1}(2i) \\
&\geq -f(2i_{0})\ (\mu_{1}-\mu_{2})\left(\{2,...,2i_{0}\}\right) + f(2i_{0})\ (\mu_{2}-\mu_{1})
\left(\{2(i_{0}+1),...,\infty\}\right)\\
\end{align*}
But\ \,$(\mu_{2}-\mu_{1})
\left(\{2(i_{0}+1),...,\infty\}\right)=-(\mu_{2}-\mu_{1})\left(\{2,...,2i_{0}\}\right)$ hence
\medskip
\begin{equation*}
\mu_{2}(f)-\mu_{1}(f)
\geq -f(2i_{0})(\mu_{1}-\mu_{2})\left(\{2,...,2i_{0}\}\right) + f(2i_{0})(\mu_{1}-\mu_{2})
\left(\{2,...,2i_{0}\}\right)\geq 0
\end{equation*}
That is why we can use theorem \ref{theo2} in this situation.
\end{rem}
We now want to apply this remark to the following probability measures on $2 \mathbb{N}-\{0\}$:
$P_{\alpha,\infty}^{0},\ P_{\alpha,h}^{\beta}$ and $P_{\alpha,h}^{\beta,+,\alpha_{1}}$ which is the
law defined in \eqref{mes2} when $\zeta_{2}\geq 0$.

First we compare $P_{\alpha,h}^{\beta}$ and $P_{\alpha,h}^{\beta,+,\alpha_{1}}$ which is in fact very easy since
\begin{align*}
&P_{\alpha,h}^{\beta,+,\alpha_{1}}\left(\tau=2\right)=P_{\alpha,h}^{\beta}\left(\tau=2\right)
\left(1+\alpha_{1}
\right)\\
&P_{\alpha,h}^{\beta,+,\alpha_{1}}\left(\tau=2r\right)=P_{\alpha,h}^{\beta}\left(\tau=2r\right)\mu_{1}
\ \,\text{for}\ r>2
\end{align*}
But $\alpha_{1}>0$ and $\mu_{1} <1$ hence
$P_{\alpha,h}^{\beta,+,\alpha_{1}}\left(\tau=2\right)>P_{\alpha,h}^{\beta}\left(\tau=2\right)$
and $P_{\alpha,h}^{\beta,+,\alpha_{1}}\left(\tau=2r\right)<P_{\alpha,h}^{\beta}\left(\tau=2r\right)$ for
$r \geq 2$. Thus remark \ref{remark} tells us that we can use theorem \ref{theo2} and define on a
 probability space
$\left(\Omega,P\right)$        a sequence of
iid random variables $\left(T^{1}_{i},T^{2}_{i}\right)_{i\geq 1}$ such that
\begin{itemize}
\item $P_{\alpha,h}^{\beta,+,\alpha_{1}}$ is the law of $T^{1}_{i}$ for every $i \geq 1$
\item $P_{\alpha,h}^{\beta}$ the law of $T^{2}_{i}$ for every $i\geq 1$
\item $P$ almost surely $T^{1}_{i}\leq T^{2}_{i}$ for every $i\geq 1$
\end{itemize}
At this point for every fixed disorder
$\zeta$ we define by recurrence another process $(T^{3}_{i})_{i\geq 1}$ with
\begin{align*}
T^{3}_{i}&=T^{2}_{i}\ \text{if}\ \zeta_{T^{3}_{1}+...+T^{3}_{i-1}+2} \geq 0\\
&=T^{1}_{i}\ \text{if}\ \zeta_{T^{3}_{1}+...+T^{3}_{i-1}+2} < 0
\end{align*}
Hence with these notations $\left(T^{2}_{i}\right)_{i\geq 1}$ is the sequence of the excursion length
of a random walk under the law $P_{\alpha,h}^{\beta}$ and $\left(T^{3}_{i}\right)_{i\geq 1}$ the one of a
random walk under the law $P_{\alpha,h}^{\beta,\zeta,\alpha_{1}}$. But by construction
$T^{3}_{i}\leq T^{2}_{i}$ for every $i\geq 1$, so for $j=\ 2\ \text {or}\ 3$ if we put
$l^{j}_{N}=\max\{s\geq 1 |T^{j}_{1}+...+T^{j}_{s}\leq N\}$ we have immediately that $P$
almost surely $l^{3}_{N}\geq l^{2}_{N}$. Thus for every $\zeta$
we have

\begin{equation*}
E_{\alpha,h}^{\beta,\zeta,\alpha_{1}}\left(\frac{l_{N}}{N}\right)=
E_{P}\left(\frac{l^{3}_{N}}{N}\right)\geq E_{P}\left(\frac{l^{2}_{N}}{N}\right)
=E_{\alpha,h}^{\beta}\left(\frac{l_{N}}{N}\right)
\end{equation*}
and integrating over $\zeta$ we obtain the left hand side of inequality \eqref{ineq}.

To finish with these inequalities we must show that the same argument allow us to compare
$\mathbb{E}\left( E_{\alpha,h}^{\beta}\left(\frac{l_{N}}{N}\right)\right)$ and
$\mathbb{E}\left( E_{\alpha,\infty}^{0}\left(\frac{l_{N}}{N}\right)\right)$. Namely we want to prove that
 remark \ref{remark} also occur. So recall

\begin{align*}
P_{\alpha,h}^{\beta}\left(\tau_{1}=2n\right)&=\left(\frac{1+\exp\left(-4hn\right)}{2}\right)\alpha^{2n}
\frac{P\left(\tau=2n\right)}{H_{\alpha,h}^{\beta}}\exp\left(\beta\right)\\
P_{\alpha,\infty}^{0}\left(\tau_{1}=2n\right)&=\frac{\alpha^{2n}P\left(\tau=2n\right)}{2 H_{\alpha,\infty}^{0}}
\end{align*}
So if we note $$L_{n}=\frac{P_{\alpha,h}^{\beta}\left(\tau_{1}=2n\right)}
{P_{\alpha,\infty}^{0}\left(\tau_{1}=2n\right)}=\left(1+\exp\left(-4hn\right)\right)
\frac{H_{\alpha,\infty}^{0}}{H_{\alpha,h}^{\beta}}\exp(\beta)$$
we immediately notice that $L_{n}$ decreases with $n$, but we have also
$$\sum_{i=1}^{\infty}P_{\alpha,h}^{\beta}\left(\tau_{1}=2i\right)=\sum_{i=1}^{\infty}
P_{\alpha,\infty}^{0}\left(\tau_{1}=2i\right)=1$$  hence necessarily there exists\ \ $i_{0}$\ \ in
$\mathbb{N}-\{0\}$ such that $P_{\alpha,h}^{\beta}\left(\tau_{1}=2i\right)
\geq P_{\alpha,\infty}^{0}\left(\tau_{1}=2i\right)$ for $i\leq i_{0}$ and
$P_{\alpha,h}^{\beta}\left(\tau_{1}=2i\right)
\leq P_{\alpha,\infty}^{0}\left(\tau_{1}=2i\right)$ for $i>i_{0}$. And the proof is complete.
\qed
\section{}
\subsection{\textit{Proof of Proposition \ref{pro:nonrand}}}
First of all, we recall a classical property which tells us that we do not
transform the free energy if we oblige the last monomer of the
chain to touch the $0$ axis. It is proved for example in a different
case in \cite{BDH} but the same technic works with our hamiltonian. So we can write
\begin{equation*}
\Phi^{0}(h,\beta)=\lim_{N\rightarrow\infty}\mathbb{E}\frac{1}{2N}\log\,E\left(\exp\left(\beta\sum_{i=1}^{2N}
\boldsymbol{1}_{\{S_{i}=0\}}-2h \sum_{i=1}^{2N}\Delta_{i}\right)\boldsymbol{1}_{\{S_{2N}=0\}}\right)\\
\end{equation*}
In the following we note
$Z_{2N,\beta,h}=E\left(\exp\left(\beta\sum_{i=1}^{2N}
\boldsymbol{1}_{\{S_{i}=0\}}-2h
\sum_{i=1}^{2N}\Delta_{i}\right)\boldsymbol{1}_{\{S_{2N}=0\}}\right)$.
\noindent Remark that $Z_{2N,\beta,h}$ can be rewrite as follow
\begin{align*}
Z_{2N,\beta,h}=&\sum_{j=1}^{N} E\left(e^{\beta
j}e^{-2h\sum_{i=1}^{2N}\Delta_{i}} \boldsymbol{1}_{\{l_{2N}=j\}}
\boldsymbol{1}_{\{S_{2N}=0\}}\right)\\
=&\sum_{j=1}^{N}\ \sumtwo{\overline l\in\mathbb N^{*j}}
{|\overline l|=N}\ \ \prod_{i=1}^{j} \left(e^{\beta j}\
V_{h,l_{j}}\right)
\end{align*}
with $V_{h,l}=P\left(\tau=2l\right)\left(e^{-4hl}+1\right)/2$.
 We aim at computing the generating function of
$Z_{2N,\beta,h}$ called $\theta_{h}(z)$
\begin{align*}
\theta_{h}(z)=&\sum_{N=1}^{\infty}Z_{2N,\beta,h} z^{2N} =\
\sum_{N=1}^{\infty}z^{2N}\sum_{j=1}^{N}e^{\beta j}\
\sumtwo{\overline l\in\mathbb N^{*j}}
{|\overline l|=N}\ \  \prod_{i=1}^{j}V_{h,l_{j}}\\
=&\sum_{j=1}^{\infty}\sum_{N=j}^{\infty}\  \sumtwo{\overline
l\in\mathbb N^{*j}} {|\overline l|=N}\ \ \prod_{i=1}^{j}
\left(e^{\beta} z^{2 l_{j}}\
V_{h,l_{j}}\right)\\=&\sum_{j=1}^{\infty}\left(\sum_{l=1}^{\infty}
e^{\beta} z^{2 l}\ V_{h,l}\right)^{j}=\
\sum_{j=1}^{\infty}\left(\sum_{l=1}^{\infty}
\frac{P(\tau=2l)}{2}\left( 1+e^{-4hl}\right) e^{\beta}
z^{2 l}\right)^{j}\\
\end{align*}
Now, recall that
\begin{equation*}
\sum_{l=1}^{\infty} P(\tau=2l) z^{2l}=1-\sqrt{1-z^{2}}
\end{equation*}
hence the computation finally gives
\begin{equation*}
\theta_{h}(z)=\sum_{j=1}^{\infty}\left(\frac{e^{\beta}}{2}\left(2-\sqrt{1-z^{2}}-\sqrt{1-z^{2}
e^{-4h}}\right)\right)^{j}
\end{equation*}
So, this serie converges when
$e^{\beta}\left(2-\sqrt{1-z^{2}}-\sqrt{1-z^{2}e^{-4h}}\right)<2$,
and if we note $R$ its convergence radius, we have
$\Phi(\beta,h)=-\log(R)$. That is why $\Phi(\beta,h)>0$ if and only
if $R<1$. So, we can say that $(h,\beta)$ is on the critical curve
if and only if for $z=1$:
$e^{\beta}\left(2-\sqrt{1-z^{2}}-\sqrt{1-z^{2}e^{-4h}}\right)=2$,
which can be write $\sqrt{1-e^{-4h}}=2\left(1-e^{-\beta}\right)$.
It gives us the critical curve equation
$$h_{c}^{0}\left(\beta\right)=\frac{1}{4}\log\left(1-4\left(1-e^{-\beta}\right)^{2}\right)$$ \qed


\bigskip

\section*{Acknowledgments}

I'm grateful to my Ph.D. supervisors Giambattista Giacomin and Roberto Fernandez
for their precious help and suggestions.


\end{document}